\documentclass{article}
\usepackage[utf8]{inputenc}
\usepackage{geometry}
\usepackage{amsmath, amsthm, amssymb}
\usepackage{amsfonts}
\usepackage{epsfig}
\usepackage{graphics}
\usepackage{subfigure}
\usepackage{graphicx}
\usepackage{color}
\usepackage{indentfirst}
\usepackage{epstopdf} %to transform eps to pdf
\usepackage{float} % for [H], to place picture in the proper position
\usepackage{mathrsfs} % For $\mathscr{L}$
\usepackage{relsize} % to use \mathlarger
\title{Divisor Functions: Train-like Structure\\ and Density Properties}
\author{Evelina Dubovski\\
Staten Island Technical High School\\
edubovski7@gmail.com}

\begin{document}

\textheight = 9in
\topmargin = -1in

\newtheorem{cor}{Corollary}
\newtheorem{df}{Definition}
\newtheorem{lm}{Lemma}
\newtheorem{rem}{Remark}
\newtheorem{thm}{Theorem}
\newtheorem{prop}{Proposition}

\newcommand{\ds}{\displaystyle}
\newcommand{\eps}{\varepsilon}
\newcommand{\noi}{\noindent}
\newcommand{\NN}{\mathbb{N}}
\newcommand{\QQ}{\mathbb{Q}}
\newcommand{\RR}{\mathbb{R}}
\newcommand{\ov}{\overline}

\newcommand{\bq}{\begin{equation}}
\newcommand{\eq}{\end{equation}}
\newcommand{\bqr}{\begin{eqnarray}}
\newcommand{\eqr}{\end{eqnarray}}

%\textheight = 9in
%\maketitle
%\date{~}
\begin{center}
{\bf \Large 
Divisor Functions: Train-like Structure and\\ Density Properties}\\[3mm]
{ Evelina Dubovski\footnote{E-mail: edubovski7@gmail.com}\\ Staten Island Technical High School}\\[2mm]
\end{center}

\begin{abstract}

We investigate the density properties of generalized divisor functions $\ds f_s(n)=\frac{\sum_{d|n}d^s}{n^s}$ and extend the analysis from the already-proven density of $s=1$ to $s\geq0$. We demonstrate that for every $s>0$, $f_s$ is locally dense, revealing the structure of $f_s$ as the union of infinitely many $trains$---specially organized collections of decreasing sequences---which we define. We analyze Wolke's conjecture that $|f_1(n)-a|<\frac{1}{n^{1-\varepsilon}}$ has infinitely many solutions and prove it for points in the range of $f_s$. We establish that $f_s$ is dense for $0<s\leq1$ but loses density for $s>1$. As a result, in the latter case the graphs experience ruptures. We extend Wolke's discovery $\ds |f_1(n)-a|<\frac{1}{n^{0.4-\varepsilon}}$ to all $0<s\leq1$. In the last section we prove that the rational complement to the range of $f_s$ is dense for all $s>0$. Thus, the range of $f_1$ and its complement form a partition of rational numbers to two dense subsets. If we treat the divisor function as a uniformly distributed random variable, then its expectation turns out to be $\zeta(s+1)$. The theoretical findings are supported by computations. Ironically, perfect
and multiperfect numbers do not exhibit any distinctive characteristics for divisor functions.
\end{abstract}
%MSC: 11N64

\section{Introduction} \label{sec1}

Perfect numbers are one of the most unique phenomena in number theory. Dating back to ancient Greece, perfect numbers first appeared in Euclid's Elements and have been the focus of many great mathematicians, including Leonardo of Pisa (Fibonacci) and Euler. But what are perfect numbers? A perfect number is a positive integer that is equal to the sum of its positive divisors, excluding the number itself.

Despite the attention that perfect numbers have received over the centuries, there is still much left to be discovered. The last known perfect number was discovered in 2018, bringing the total of identified perfect numbers to 51. We do not even know whether there exist odd perfect numbers or whether their cardinality is infinite. An analysis of perfect numbers has led us to our investigation of the divisor function and its properties. The intended audience is researchers in multiplicative number theory and the broad community of math lovers.

The function $\sigma(n)=\sum\limits_{d|n} d$ is called the divisor function of whole number $n$. It is equal to the sum of all divisors $d$ (factors) of $n$, including 1 and $n$.
For example, $\sigma(10)=18$. Clearly, $\sigma(p)=1+p$ if $p$ is prime. 
%If a number is equal to the sum of its divisors (without the original number), then it is called {\it perfect}.
We consider the generalization of perfect numbers in terms of the divisor function. A number is perfect if $\sigma(n)=2n$. Since $\sigma(6) = 12$, then 6 is perfect along with 28, 496, etc. Scientists also consider multiperfect numbers, for which divisor function $\sigma$ takes integer values called abundancies. If $\sigma(n)= kn$, then $n$ is multiperfect of abundancy $k$. A recent work \cite{2008multiperfect} states that no $k$-perfect
odd numbers are known for any $k\geq2$. 

The motivation for this research was to analyze whether perfect and multiperfect numbers are really special, leading, as a research tool, to the investigation of the properties of divisor functions presented in this article. 
%was the question whether perfect numbers are really specific. This question led towards the investigation of the properties of divisor functions presented in this article.} 
The theory of the divisor function $\sigma(n)$ and its relation to other challenging problems can be found, e.g., in books \cite{2008Hardy} and \cite{1999Landau}. Interesting recent results on the sums involving the divisor function are in \cite{2015Pollack, 2021Ma-Sun}. 
%\cite{2008Hardy, 2021Ma-Sun, 2009Akbary, 2015Pollack}.

In 1941, Cramer \cite{1941Cramer} proved that the function $\displaystyle  f(n)=\frac{\sigma(n)}{n}$ is dense on $[1,\infty)$, i.e., for every $a\in[1,\infty)$ and any positive small $\varepsilon>0$, there is at least one $n$ such that $|f(n)-a|<\varepsilon$. Since $\varepsilon$ is arbitrary, then “at least one” actually implies “infinitely many.”

In 1977, Wolke \cite{1977Wolke} strengthened 
this result as follows: for every $a\geq1$ and every $\varepsilon>0$, the inequality 
\bq
\left|f(n) -a \right|\leq \frac1{n^{0.4-\varepsilon}} \label{Wolke}
\eq
has infinitely many solutions in whole numbers $n$. This is further evidence that the function $\displaystyle f(n)=\frac{\sigma(n)}{n}$ is dense in $[1,\infty)$. In fact, the above inequality implies that in any small neighborhood of any real $a\geq1$ we can find infinitely many values of $f(n)$. Moreover, in his paper %\cite{1977Wolke} 
D. Wolke conjectured that constant $0.4$ can be replaced by 1.
%We prove his conjecture for $a$ from the range of $f$.\\

Our goal is to examine the generalization of the divisor function,\linebreak 
$\sigma_s(n)=\sum_{d|n} d^s$ for $s\geq0$.
If $s=1$, we arrive at 
the usual divisor function $\sigma(n)$. If $s=0$, we obtain $\sigma_0(n)$, which is just the number of divisors of $n$. 

We investigate the function
$f_s(n)=\displaystyle\frac{\sigma_s(n)}{n^s}$ for different non-negative values of $s$ to determine whether their ranges are dense (like $f_1$) or not. We explore the structure of 
functions $f_s$ and show that their range is the union of trains, specially organized linked decreasing sequences, that are defined in this paper. This observation leads to the proof of Wolke's conjecture for the points from $R(f_s)$, the range of $f_s$. 
That is, we show that the extended Wolke's conjecture is valid for any $a\in R(f_s)$. Particularly, for $s=1$, the constant 0.4 in inequality (\ref{Wolke}) can be replaced by 1.

For all reals, not only from the range, 
we rely on Cramer's approach \cite{1941Cramer} and demonstrate that for $0<s\leq 1$, the function $f_s$ is dense and chaotic, but for $s>1$, $f_s$ loses density and ruptures appear. 

%Based on this analysis, we prove the following generalization of Wolke's inequality:
%\bq
%\left|f_s(n) -a \right|\leq \frac1{n^{0.4s-\varepsilon}}, \quad 0<s\leq1. \label{Wolke2}
%\eq
 
By extending Wolke's finding (\ref{Wolke}) to  $0<s\leq 1$, we obtain the following quantitative measure of the density strength for all $a\geq1$:
\bq
\left|f_s(n) -a \right|\leq \frac1{n^{0.4s-\varepsilon}}.  \label{Wolke2}
\eq

We end the paper by proving that the rational complement to the range is also dense in $[1,\infty)$ for positive integer $s$. 
The presented analysis is supported by computations.\\

%%%%%%%%%%%%%%%%%%%%%%%%%%%%%%%%%%%%%%%%%%%%%%%%%%%%%%%%%%%%%%%%%%%%%%%%%%%%%%%%

\section{Local density} \label{sec2}
%The goal of this section is to show, that $f_1$ is actually denser with the density strength equal to $1$ and estimate the density strength for different $s$.
The most essential property of the divisor function is its multiplicativity for the products of relatively prime numbers. If $a$ and $b$ are relatively prime and $c_i$ are the factors of $a$ whereas $d_i$ are the factors of $b$, then
\[
\sigma_s(ab)= \sum_{i,j} c_i^sd_j^s = \sum_i c_i^s\cdot \sum_j d_j^s = \sigma_s(a)\cdot \sigma_s(b).
\]
Particularly, if we consider the prime factorization $n=\prod_{i=1}^m p_i^{k_i}$, then $\sigma_s(np)=\sigma_s(n)\cdot \sigma_s(p)$ if prime number $p$ does not coincide with any of $p_i$, i.e, $p\not\in P_n=\{p_1,p_2,\ldots,p_n\}$, the set of all prime factors of $n$. Similar results hold for functions $f_s$, and we have $f_s(np)=f_s(n)\cdot f_s(p)$. So,
\begin{equation}
f_s(n)=f_s\left(\prod_{i=1}^m p_i^{k_i}\right)=\prod_{i=1}^m f_s(p_i^{k_i}). \label{1_}
\end{equation} 

The computations shown below indicate that the graphs of $f_s(n)$ are formed by an infinite union of subgraphs, each of which is a sequence that
decreases as $n\to\infty$. Motivated by these computational results, we consider prime $p\not\in P_n$ and generate the decreasing sequence $\{f_s(np)\}_{p\not\in P_n}$.
%$x_{s}[n]$,   
%\[
%x_{s}[n]=\left\{ f_s(np):\ p\not| n\right\}.
%\]
In the well-investigated case $s=1$, we obtain for $n=2$ the sequence  
$\{f_{1}(2p)\}=\left\{2, 9/5, 12/7, 18/11,\ldots \right\}$, where the corresponding values are $p=3, 5, 7, 11,\ldots$.
For $n=6$, we obtain the sequence  
$\{f_{1}(6p)\}=\left\{12/5, 16/7, 24/11, \ldots \right\}$. Here $p=5, 7, 11,\ldots$. Clearly, these sample sequences are decreasing.

\begin{lm} \label{lm1}
 For $s>0$, the sequences $\{f_s(np)\}$ are decreasing for all sufficiently large $p$ and \linebreak $\lim_{p\to\infty} f_s(np)=f_s(n)$.
 \end{lm}
\noi{\bf Proof.} Since $p\to\infty$, we can assume that $p\not \in P_n$. Then 
\begin{equation}
f_s(np)=f_s(n) f_s(p)= f_s(n)\cdot \left (1+\frac1{p^s}\right) \to f_s(n)\quad {\rm as}\ p\to\infty.
\label{2}
\end{equation}
This proves Lemma \ref{lm1}. \hfill $\qed$\\

We introduce the following definition to conveniently describe the structure of the range of $f_s$.
\begin{df} \label{train}
We call {\rm train} an ordered collection of linked sequences such that\\
{\rm(1)} any term of the sequence is greater that any term of the preceding sequence;\\
{\rm(2)} the infimum of any sequence is equal to the supremum of the preceding sequence (i.e., the sequences are {\rm linked}).\\
The above sequences, forming the train, are called {\rm cars}. 
\end{df}
This terminology comes from real trains in which the passengers of any car are ahead of all passengers in the following car.

\begin{rem} If $p\in P_n$, then sequence $\{f_s(np)\}$ is not necessarily decreasing for small $p$. For example, if $n=2$ and $p=2$ then $\{f_1(np)\} = \{1.75, 2, 1.8,\ldots\}$.
\end{rem}

\begin{thm} \label{thm1}
{\rm {\bf (Local density)} Let $s>0$. If $r=f_s(n)$ is a number from the range $R(f_s)$ of function $f_s$, then for every $\varepsilon>0$ there are infinitely many values of $m$ such that $f_s(m)\in[f(n),f(n)+\varepsilon]$.}
\end{thm}
\noindent The proof of Theorem \ref{thm1} follows directly from Lemma \ref{lm1}.\\

\begin{figure}[H] 
	\begin{center}
	\includegraphics[width=6in]{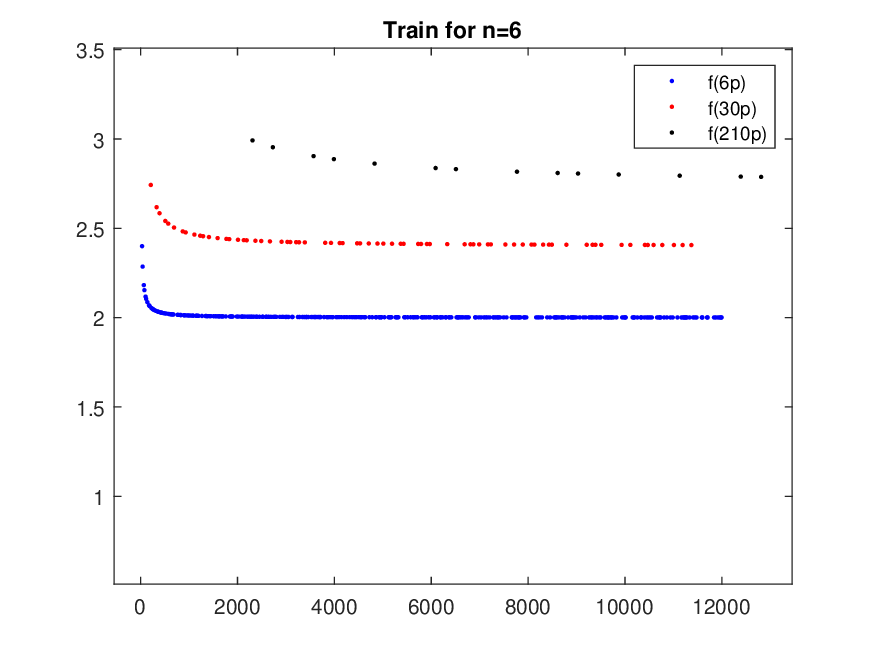}
				\caption{
					\label{Fig1} Three first sequences (“cars") of the train for $n=6$, $s=1$.}
    We begin with the blue sequence $f(6p)$, in which the first value of $p$ is 5, so $f(6\times 5)$ is the highest value in this “car." The decreasing sequence $f(6p)\to f(6)$ as $p\to \infty$. Any term of the next “car" is greater than any term of the previous sequence. As shown, $f(30p)=f(6\times 5\times p)$ decreases to $f(30)$, the highest term of $f(6p)$. Once again, any term of a “car" is greater than any term of the previous “car." Just as the cars of a train are connected but do not overlap, so are these sequences.
    %}	\label{Train}
	\end{center}
\end{figure}

As such, it can be said that the range of $f_s$ is “dense from above.” Also, the function is locally dense: in any neighborhood of any point $r\in R(f_s)$, there are infinitely many points from range $R(f_s)$.
Figure \ref{Fig1} clarifies the statement of Theorem \ref{thm1}.

According to the proof of Lemma \ref{lm1}, for a given $n=p_1^{k_1}p_2^{k_2}\cdot\ldots\cdot p_m^{k_m}$, if $q_1$ is the smallest prime outside the set $P_n$, then the sequence $\{f_s(np)\}$ attains its maximum at $p=q_1$, which is $\displaystyle f_s(n)\cdot(1+\frac1{q_1^s})$. Similarly, the sequence $\{f_s(nq_1p)\}$ attains its maximum at $p=q_2$, which is the lowest prime number outside $P_{nq_1}$. By continuing, we arrive at the sequence
$\{f_s(nq_1q_2\cdot\ldots\cdot q_jp)\}_{p\not\in P_n}$. 
Here, the primes $q_i\not\in P_n$ and form the ordered sequence outside of $P_n$, $q_1<q_2<\ldots q_j$. The graphs for $f_s$ are described by the following assertion:
\begin{prop}\label{prop1}
For every prime $p>q_{j+1}$, $p\not\in P_n$,
\[
 f_s(np)<f_s(nq_1p)<\ldots<f_s(nq_1q_2\cdot\ldots\cdot q_jp).
\]
Thus, every term from sequence $\{f_s(nq_1q_2\cdot\ldots\cdot q_jq_{j+1}p)\}$ with variable primes $p$ is greater than any term from sequence $\{f_s(nq_1q_2\cdot\ldots\cdot q_{j}p)\}$. 
\end{prop}
\noi{\bf Proof.} In view of (\ref{2}), 
\begin{eqnarray*}
&& f_s(nq_1q_2\cdot\ldots\cdot q_jp)=
f_s(nq_1q_2\cdot\ldots\cdot q_j)f_s(p)\\
&&\phantom{WW} < f_s(nq_1q_2\cdot\ldots\cdot q_j)f_s(q_{j+1}) =
f_s(nq_1q_2\cdot\ldots\cdot q_jq_{j+1})\\
&&\phantom{WW}< f_s(nq_1q_2\cdot\ldots\cdot q_jq_{j+1}p).
\end{eqnarray*}
Therefore, 
\[
\max_p f_s(nq_1\ldots q_jp)\le f_s(nq_1\ldots q_jq_{j+1})<f_s(nq_1\ldots q_jq_{j+1}p).
\]

\noindent This proves Proposition
\ref{prop1} \hfill \qed \\

Proposition \ref{prop1} explains the visible structure of the graphs for functions $f_s$ (see Figures \ref{Fig1}, 2, and 3). Each of the sequences can be considered a “car" on the train, or collection, of sequences $\{f_s(nq_1q_2\cdot\ldots\cdot q_jp)\}$.

%%%%%%%%%%%%%%%%%%%%%%%%%%%%%%%%%%%%%%%%%%%%%%%%%%%%%%%%%%%%%%%%%%%%%%%%%%%%%%%%%

Relying on the above analysis, we now strengthen Wolke's estimate (\ref{Wolke}) for the points from the range $R(f_s)$.
\begin{thm} \label{thm_Wolke}
Let $s>0$. For every $a\in R(f_s)$ there are infinitely many solutions to inequality
\bq \label{Wolke_range}
|f_s(n)-a|\leq \frac{c(a)}{n^s}.%\frac{af_s^{-1}(a)}{n^s}.
\eq
\end{thm}
\noindent {\bf Proof.} Let $a=f_s(m)$. Then for every $p\not\in {P}_m$ we have $\ds f_s(mp)=a\cdot \left(1+\frac1{p^s}\right)$. 
Hence, for all $n=mp$ and $p\not\in {P}_m$, we obtain
$\ds f_s(n)-a=\frac{a[f_s^{-1}(a)]^s}{n^s}$. This proves Theorem~\ref{thm_Wolke}. \hfill $\qed$ \\
\begin{cor} Let the conditions of Theorem \ref{thm_Wolke} hold. Then for every $\eps>0$, the inequality
\bq \label{Wolke_range2}
|f_s(n)-a|\leq \frac1{n^{s-\eps}},\quad s>0,%\frac{af_s^{-1}(a)}{n^s}.
\eq
has infinitely many solutions. 
\end{cor}
\noindent {\bf Proof.} The shift by $\eps$ compensates for the constant $c(a)$ and allows us to replace it by 1. \hfill $\qed$

%%%%%%%%%%%%%%%%%%%%%%%%%%%%%%%%%%%%%%%%%%%%%%%%%%%%%%%%%%%%%%%%%
% SECTION 3
\section{Global density} \label{sec3}
Let $n=p^k$ for prime $p$. Then 
\bq
f_s(p^k)=\frac{1+p^s+\ldots+p^{ks}}{p^{ks}}=\frac{p^{(k+1)s}-1}{(p^s-1)p^{ks}}=\frac{p^s-p^{-ks}}{p^s-1}=\frac{p^s}{p^s-1} - \frac{p^{-ks}}{p^s-1}. \label{3-1}
\eq
In general, for $n=\prod\limits_{i=1}^{m} p_i^{k_i}$, we use multiplicative property (\ref{1_}) and obtain
\begin{equation}
f_s(n)=f_s\left(\prod_{i=1}^{m} p_i^{k_i}\right)=\prod_{i=1}^{m} f_s\left(p_i^{k_i}\right)
=\prod_{i=1}^{m} \left(\frac{p_i^s}{p_i^s-1} - \frac{p_i^{-k_is}}{p_i^s-1}\right). \label{3}
\end{equation}
The second term $p_i^{-k_is}\cdot(p_i^s-1)^{-1}$ becomes arbitrarily small for big $p_i$. 
The first, principal, term in (\ref{3}) yields
\begin{equation}
f_s(n)\approx \prod_{i=1}^{m} \frac{p_i^s}{p_i^s-1}=Q_m. \label{4}
\end{equation}
Let's investigate how closely $Q_m$ approximates $f_s(n)$.\\
It is known \cite{1999Landau} that 
\begin{equation}
\lim_{m\to\infty} \prod_{i=1}^{m} \left( 1+\frac{1}{p_i-1}\right) = \infty,  \label{5}
\end{equation}
where $p_i$ are all prime numbers.

%% LEMMA 2
\begin{lm}\label{lm2}
For $0<s\leq 1$,

\begin{equation}
\prod_{i=1}^\infty \left( 1+\frac{1}{p_i^s-1}\right) = \infty  \label{5'}
\end{equation}
 whereas for $s>1$ this infinite product is finite
\begin{equation}
\prod_{i=1}^\infty \left( 1+\frac{1}{p_i^s-1}\right) < \infty. \label{5a}
\end{equation}
\end{lm}
\noindent{\bf Proof.}
For $0<s\leq 1$ we obtain
\[
1+\frac{1}{p_i^s-1}\geq 1+\frac{1}{p_i-1}=\frac{p_i}{p_i-1}.
\]
In view of (\ref{5}) we obtain (\ref{5'}): 
\begin{equation}
\prod_{i=1}^{\infty} \left(1+\frac{1}{p^s_i-1}\right) \geq \prod_{i=1}^{\infty} \left(1+\frac{1}{p_i-1}\right) = \infty.\label{6}
\end{equation}

\noindent For $s>1$ we have 
\bq \label{9_}
\prod_{i=1}^{\infty} \left( 1+\frac{1}{p_i^s-1}\right)=\exp\left(
\sum_{i=1}^\infty \ln\left(1+\frac{1}{p_i^s-1}\right)\right).
\eq
Since $\ln(1+x)\leq x$, then we consider
\bq \label{10}
\sum_{i=1}^\infty \ln\left(1+\frac{1}{p_i^s-1}\right)\leq 
\sum_{i=1}^\infty \frac{1}{p_i^s-1}<1+\sum_{i=2}^\infty \frac{1}{p_i^s-1}.
\eq
Since $p_i^s\geq (p_{i-1}+1)^s>p_{i-1}^s+1$, then $\ds \frac1{p_{i}^s-1}<\frac1{p_{i-1}^s}$, and for (\ref{10}) we obtain the upper bound
\[
\sum_{i=1}^\infty \frac{1}{p_i^s-1} < 1+\sum_{i=2}^\infty \frac1{p_{i-1}^s} = 1+\sum_{i=1}^\infty \frac1{p_{i}^s} < \sum_{i=1}^\infty \frac1{i^s}<\infty.
\]
Combining this result with (\ref{9_}), (\ref{10}), we derive (\ref{5a}).\hfill $\qed$\\
%Another way to
%prove (\ref{5a}) is to apply the Euler product formula
%\[
%\prod_{i=1}^\infty \left(1+\frac1{p_i^s-1}\right) = \zeta(s) 
%\]
%and use the boundedness of zeta-function $\zeta(s)=\sum_{i=1}^\infty \frac1{i^s}$ at $s>1$.\hfill $\qed$\\

%%%%%%%%%%%%% Lemma 3

In Section \ref{sec2}, we showed that the points from range $R(f_s)$ can be approximated from above. Now we consider arbitrary real numbers $a$ and show that $a>1$ can be approximated from below.

\begin{lm} \label{lm3} 
Let $0<s\leq1$. Then for any real number $a>1$ and any $\varepsilon>0$, there exist primes $q_1<q_2<\ldots<q_m$ such that 
$a-\eps<Q_m  <a$, where
\[
Q_m=\prod_{i=1}^{m}\left( 1+ \frac{1}{q_i^s-1}\right).
\]
\end{lm}
\noi{\bf Proof.} Let $q_1$ be the smallest prime such that $\displaystyle 1+\frac{1}{q^s_1-1}\leq a$. Since the infinite product  (\ref{5}) diverges and
$\ds 1+\frac{1}{p^s_i-1}\to 1$ as $p_i\to\infty$, we can choose additional $j_1$ consecutive bigger primes $q_i$, $1\leq i\leq j_1$ such that 
the corresponding product 
\[
Q_{j_1}=\prod\limits_{i=1}^{j_1} \left(1+\frac1{q_i^s-1}\right)
\]
is still less than $a$. Let $q_{j_1}$ be the last $q$ which can be chosen this way. Then $Q_{j_1}<a$ but 
\begin{equation}
Q_{j_1}\cdot \left(1+\frac1{\hat{q}^s-1}\right)\geq a \label{6a}
\end{equation}
where $\hat{q}$ is the nearest prime number greater than $q_{j_1}$.\\
Then we skip several prime numbers and select the smallest prime $q_{j_1+1}$ such that 
\[
Q_{j_1+1}=Q_{j_1}\cdot \left(1+\frac1{q_{j_1+1}^s-1}\right)\leq a.
\]
Once again, we select the next finite set of consecutive primes 
$q_{j_1+1}$, $q_{j_1+2}$, $\ldots$, $q_{j_2}$ such that $Q_{j_2}\leq a$ but $\ds Q_{j_2}\cdot \left(1+\frac1{\hat{q}^s-1}\right)\geq a$. As before, $\hat{q}$ is the nearest prime number greater than $q_{j_2}$. 
Then we repeat the procedure, finding $q_{j_3}$, $q_{j_4}$, $\ldots$, $q_{j_k}$ such that the corresponding 
\[
Q_{j_3}< Q_{j_4} < \ldots < Q_{j_k}\leq a
\]
and, eventually, $Q_{j_k} \geq a-\eps$. \\
Thus, $m=j_k$, proving Lemma \ref{lm3}. So, we can arrive arbitrarily close to $a$ from below. \hfill  $\qed$\\

%
%The next question is whether $Q_j<a$ is sufficiently close to $a$ from below. 
%To arrange the proper proximity, we use another method for the choice of $q_i$ for $i>j$. By Bertrand's postulate, there is a prime between $x$ and $2x-2$ for any $x>3.5$.

%Let us pick up $b_j\in\RR$ such that $Q_j\left(1+\frac1{b_j}\right)=a$. In view of (\ref{6a}), it follows that $b_j\geq q_{j+1}^s-1$.\\

%Obviously, if $a>\frac{3^s}{3^s-1}$, we can use $q_j=q_1=3$. 
%If $a\leq \frac{3^s}{3^s-1}$, we have $q_j>3$. In both cases $b_j\geq 5^s-1$ since the next prime greater than $q_{j}$ is greater than 3. Then there is a prime $q_{j+1}$ such that $b_j$..... 

%then if we let $q_1=p_m$, $q_2=p_{m+1}$, $q_3=p_{m+2}$, and so on, we eventually approach arbitrarily close to $a$ from below. 
%cause $Q_i$ to exceed $a$ for some $i$. But we continue this rule for choice of $q_i$ only as long as $Q_i<a$ and then choos

%%%%%%%%%%%%% Lemma 4
The following lemma is needed in order to analyze the negligibility of the small term 
$\ds \frac{p_i^{-k_is}}{p_i^s-1}$ in the infinite product (\ref{3}). 
\begin{lm}\label{lm4}
Let $s>0$. For any set of $m$ primes $q_1<q_2<\ldots q_m$ and any $\varepsilon>0$ there exist exponents 
$n_1,n_2,\ldots, n_m$ such that 
\begin{equation}
Q_m-\varepsilon<f_s\left(q_1^{sk_1}q_2^{sk_2}\times\ldots\times q_m^{sk_m}\right)<Q_m \label{7}
\end{equation}
whenever $k_i>n_i$.
\end{lm}

\noi{\bf Proof.} Following (\ref{3-1}), let us consider $\ds f_s(q^{sx})=\frac{q^{s(x+1)}-1}{q^{sx}(q^s-1)}$. We differentiate in $x$ and obtain $\ds \frac{d}{dx}f_s(q^{sx})=\frac{s\ln q}{(q^s-1)q^{sx}}>0$.
Thus, $f_s$ increases with $x$. Since $\ds \lim_{x\to\infty} f_s(q^{sx})=\frac{q^s}{q^s-1}$, then for any $\varepsilon>0$, there is an $n$ such that 
\[
\frac{q^s}{q^s-1}-\eps< f_s(q^{sk})<\frac{q^s}{q^s-1}
\] 
 whenever $k>n$. It follows that, for any two primes $q_1$ and $q_2$ and any $\varepsilon_1>0$ and $\varepsilon_2>0$, there are numbers $n_1$ and $n_2$ such that for all $k_1>n_1$ and $k_2>n_2$ the following inequalities hold:
\begin{equation}
\frac{q_i^s}{q_i^s-1}-\varepsilon_i <f_s\left(q_i^{sk_i}\right) < \frac{q_i^s}{q_i^s-1}. \label{8_}
\end{equation}
Using the multiplicative property, $f_s(q_1^{sk_1} q_2^{sk_2})=f_s(q_1^{sk_1})f_s(q_2^{sk_2})$, we obtain by multiplication the inequalities 
(\ref{8_}) for $i=1$ and $i=2$:
\begin{eqnarray}
Q_2-\left(\frac{\varepsilon_1 q^s_2}{q^s_2-1} +\varepsilon_2 Q_1-\varepsilon_1\varepsilon_2\right)<f_s(q_1^{sk_1} q_2^{sk_2})<Q_2. \label{99}
\end{eqnarray}
If we choose $\ds\varepsilon_1<\varepsilon\frac{q^s_2-1}{2q^s_2}$ and $\ds\varepsilon_2<\frac{\varepsilon}{2Q_1}$, then from (\ref{99}) we have
\bq
Q_2-\varepsilon<f_s(q_1^{sk_1} q_2^{sk_2})<Q_2. \label{15}
\eq
Then we multiply (\ref{15}) by (\ref{8_}) at $i=3$, then by (\ref{8_}) at $i=4$, and so on, extending this argument to the complete set of $m$ primes by induction. \hfill $\qed$\\  

%%%%%%%%%%%%%%%%% THEOREM 3
\begin{thm}\label{thm2}
Let $0<s\leq1$. Then there exist infinitely many integers $n$ such that $f_s(n)$ differs from $a$ by the amount less than arbitrary $\varepsilon>0$ where $a$ is any real constant, $a\geq1$.
\end{thm}
\noindent{\bf Proof.} By Lemma \ref{lm3} we take $m$ primes $q_1,q_2,\ldots, q_m$ such that $Q_m$ differs from $a$ by an amount less than $\frac12\varepsilon$. Lemma \ref{lm4} allows us to choose the 
exponents of these primes such that $f_s\left(q_1^{sk_1}q_2^{sk_2}\ldots q_m^{sk_m}\right)$ differs from $Q_m$ by less than $\frac12\varepsilon$. 
Finally, $\left|f_s\left(q_1^{sk_1}q_2^{sk_2}\ldots q_m^{sk_m}\right) - a\right|<\varepsilon$, proving Theorem \ref{thm2}.\hfill $\qed$\\

The computations confirm the structure of the range $R(f_s)$ as the union of trains. Figure \ref{Fig1} showed the train $T(6)$. Trains are also observable in Figures 2 and 3 below. 

In Figure 2 ($s=1$), the trains look dense, fitting Theorem \ref{thm2}. In Figure 3 the computations verify that density fails for $s>1$ as proved in Theorem \ref{thm3A} below. 
%Thus, the range of $f_s(n)$ is dense for $0<s\leq1$ and has ruptures otherwise.

%The presented computational figures justify these conclusions. 

~\hspace{-2cm}    \vspace{-1cm} \includegraphics[height=8cm]{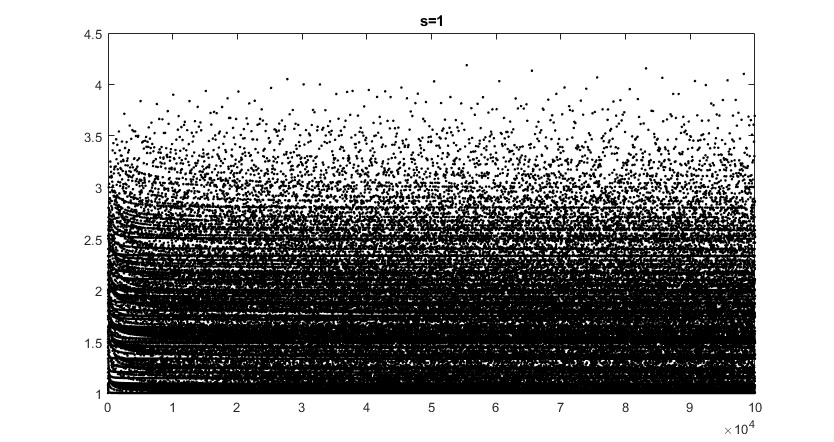}
  ~\hspace{-1cm} 
    \includegraphics[height=8cm]{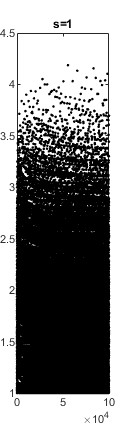}\\[2mm]
\begin{center} Figure 2: Dense behavior of $f_s(n)$ at $s=1$. A “chaos" structured by trains (see Fig. \ref{Fig1}).\\ 
The compressed graph on the right makes density more visible.
\end{center}

%%%%%%%%%%%%%%%%%%%%%%%%%%%%%%%%%%%%%
% THEOREM 4
%The following result proves that for $s>1$ the density property fails.
\begin{thm} \label{thm3A}
In order for the functions $f_s$ to lose the density property and have ruptures, it is necessary and sufficient that $s>1$.
\end{thm}
\noi{\bf Proof.} The necessity follows from 
the key point of the proof of Theorem \ref{thm2}, the divergence of the infinite product (\ref{6}). For $s>1$, the divergence fails.\\
To prove sufficiency, let us find the upper bound for range $R(f_s)$. From (\ref{3-1}),
\bqr
f_s(n)=\prod\limits_{i=1}^m\frac{1-p_i^{-(k+1)s}}{1-p_i^{-s}}%=\prod\limits_{i=1}^m\left(1+\frac{1-p_i^{-ks}}{p_i^s-1}\right) 
< \prod\limits_{i=1}^\infty \frac1{1-p_i^{-s}} = \sum_{n=1}^\infty \frac1{n^s}.\nonumber 
\eqr
%\left(1+\frac{1}{p_i^s-1}\right)\leq \exp\left(\sum_{i=1}^\infty \frac1{p_i^s-1}\right) < \exp\left(\sum\limits_{n=1}^\infty \frac1{n^s}\right).
The last step comes from Euler's identity.  
%inequality 
%\[
%(n+1)^s - n^s - 1 > 0,\quad s>1.
%\]
Finally, for every $n\in\NN$ and $s>1$,
\bq
f_s(n) < \zeta(s)<\infty. \label{upperbound} \hfill \qed
\eq
Particularly, for $s=2$ we have $f_2(n)<\frac{\pi^2}{6}$. The computations (Figure 3) below perfectly fit this estimate and show that density fails for $s>1$. 

~\hspace{-2.1cm}    \vspace{-1cm} \includegraphics[height=8cm]{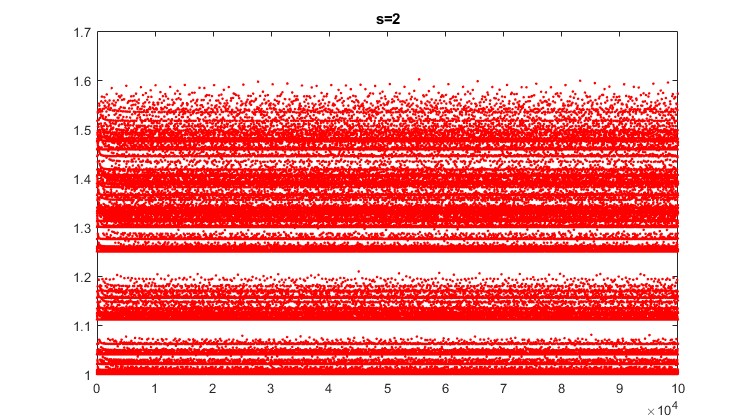}
  ~\hspace{-1cm} 
    \includegraphics[height=8cm]{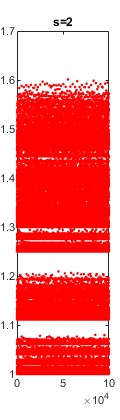}\\[3mm]
\begin{center} Figure 3: Non-dense range of $f_s$ at $s=2$. $f_2(n)<\frac{\pi^2}{6}$. No density and the appearance of ruptures for $s>1$. Trains are also observable.
The compressed graph of $f_2$ on the right makes the ruptures more visible.
\end{center}

%%%%%%%%%%%%%%%%%%%%%%%%%%%%%%%%%%%%%%%%%%%%%%%%%%%%%%%%
% THEOREM 5
The following Theorem strengthens the density result of Theorem \ref{thm2} and provides the quantitative measure of the density. Here, we extend inequality (\ref{Wolke2}) from $s=1$ to all $0<s\leq1$.
\begin{thm}\label{thm4}
Let $0<s\leq 1$. Then for every $a\geq1$ and $\eps>0$ the inequality
\bq
|f_s(n)-a|\leq \frac1{n^{0.4s-\eps}} \label{wolke3}
\eq
has infinitely many solutions on $[1,\infty)$. 
\end{thm}
\noindent{\bf Proof.} We follow Wolke's idea \cite{1977Wolke}
and use the following notations: $\ov{p}$ is the greatest prime less than prime number $p$ whereas $\hat{p}$ is the smallest prime greater than $p$. So, $\ov{p}<p<\hat{p}$. We choose any positive $y$ and find prime $p_0>y$ such that
\bq
\prod_{y<p\leq p_0} \left(1+\frac1{p^s}\right)\leq a\quad {\rm and}\quad \left(1+\frac1{\hat{p}^s}\right)\cdot\prod_{y<p\leq p_0} \left(1+\frac1{p^s}\right)>a \label{w1}
\eq
Let us consider the product of primes
\[
n_0=\prod_{y<p \leq p_0} p. 
\]
Similar to the previous constructions, let $p_1$ be the smallest prime such that
$\ds f_s(n_0)\left(1+\frac1{p_1^s}\right)<a$. 
We define step-by-step the sequence $\{n_k\}$ as follows: 
\bq
n_k=n_0 p_1p_2\cdot\ldots\cdot p_k,\quad f_s(n_k)<a,\quad {\rm but}\  f_s(n_{k-1})\cdot\left(1+\frac1{\ov{p}_k^s}\right)>a. \label{w3}
\eq
Let $p_{k+1}$ be the smallest prime number such that
\[
f_s(n_k)\cdot \left(1+\frac1{p^s_{k+1}}\right) \leq a.
\]
Thus, $p_{k+1}>p_k$ and $n_{k+1}=n_k p_{k+1}$. Certainly, if $f_s(n_{k+1})=a$, then the construction of the sequence $\{n_k\}$ ends. Further, we consider infinite sequence $\{n_k\}$.\\
Due to the above construction of $\{n_k\}$, it is easy to observe 
\[
f_s(n_{k+1})=f_s(n_{k-1})\left(1+\frac1{p^s_k}\right)\cdot \left(1+\frac1{p^s_{k+1}}\right) \leq a\leq f_s(n_{k-1})\left(1+\frac1{\ov{p}^s_k}\right).
\]
Then
\[
\ds 1+\frac1{p^s_{k+1}} \leq \frac{\ds 1+\frac1{\ov{p}^s_k}}{\ds 1+\frac1{p^s_k}}
=\left(1+\frac1{\ov{p}^s_k}\right) \sum_{i=0}^\infty (-1)^{i} \frac1{p^{is}_k}
\leq 
\left(1+\frac1{\ov{p}^s_k}\right) \left(1- \frac1{p^{s}_k} + \frac1{p^{2s}_k}\right).
\]
%We used here $\ds \frac1{1+x}\leq \frac43(1-x)$ for $x\leq\frac12$.
Since $\ov{p}_k>\frac12 p_k$, then we obtain
\bq
p_{k+1}^s\geq \frac{\ov{p}_k^s p_k^s}{p_k^s-\ov{p}_k^s}\geq   \frac1{2^s}\cdot\frac{p_k^{2s}}{p_k^s-\ov{p}_k^s}. \label{w4}
\eq
We introduce $\Delta_k$ such that $f_s(n_k)\cdot \left(1+\frac1{\Delta_k}\right) = a$. Then from the inequality
\[
f_s(n_k)\left(1+\frac1{\ov{p}_{k+1}^s}\right)>a
\]
we obtain
\[
1+\frac1{\Delta_k} < 1+\frac1{\ov{p}_{k+1}^s}.
\]
Taking (\ref{w4}) into account, we get
\bq
\Delta_k > \ov{p}^s_{k+1}\geq \frac1{2^s} p^s_{k+1}\geq \frac1{4^s}\cdot \frac{p_k^{2s}}{p_k^s-\ov{p}_k^s}.
\nonumber
\eq
Next,
\bq
\left|f_s(n_k) - a\right| = \left|f_s(n_k) - f_s(n_k)\cdot \left(1+\frac1{\Delta_k}\right)\right|
\leq \frac{f_s(n_k)}{\Delta_k}\leq 4^s a\cdot \frac{p^s_k - \ov{p}^s_k}{p^{2s}_k}. \label{w5}
\eq
As Montgomery \cite{1971Mont} showed, for any $\delta>0$ and sufficiently big $x>x_0(\eps)$, there is a prime in the interval $[x,x+x^{0.6+\delta}]$. Then 
\[
p^s-\ov{p}^s\leq p^{0.6s+\delta s}% - replaced \leq 2^{-s} p^{0.6s+\delta s} thanks to the choice of delta.
\]
From (\ref{w4}), for sufficiently big $k\geq k_0(\delta)$, we have $p^s_{k+1}\geq p_k^{(1.4-\delta) s}$. 
Then,
\bq
p_{k+2}^s\geq p_{k+1}^{(1.4-\delta) s}\geq p_k^{(1.4-\delta)^2 s^2}, \label{w5a}
\eq
and we have 
\[
p_{k+m}^s\geq p_k^{(1.4-\delta)^m s^m}.
\]
Therefore,
\[
p_{k-m}^s\leq p_{k}^{(1.4s-\delta s)^{-m}}. 
\]
This yields 
\bq \label{w6}
n_k\leq n_{k_0} p_k^{(0.4s-\delta  s)/(1.4s - \delta s)}.
\eq
In addition, using (\ref{w5}), (\ref{w5a}), 
\[
|f_s(n_k)-a|\leq C \frac{p_k^{0.6s+\delta s}}{p^{2s}_k} =C p_k^{-1.4s+\delta s}.
\]
Then we use (\ref{w5}) and (\ref{w6}) to obtain
\[
|f_s(n_k)-a|\leq C p_k^{-1.4s+\delta s}\leq n_k^{-0.4s(1+\eps)}\leq \frac1{n_k^{0.4s-\eps}}.
\]
Since there are infinitely many $n_k$, inequality (\ref{wolke3}) has been proved.  
Thus, we can state that Theorem \ref{thm4} generalizes Wolke's discovery. The shadowed area below demonstrates the essence of Theorem \ref{thm4} and contains infinitely many solutions to inequality (\ref{wolke3}).\\

%~\hspace{-2.1cm}    \vspace{-1cm} 
\includegraphics [width=12cm, height=6cm]{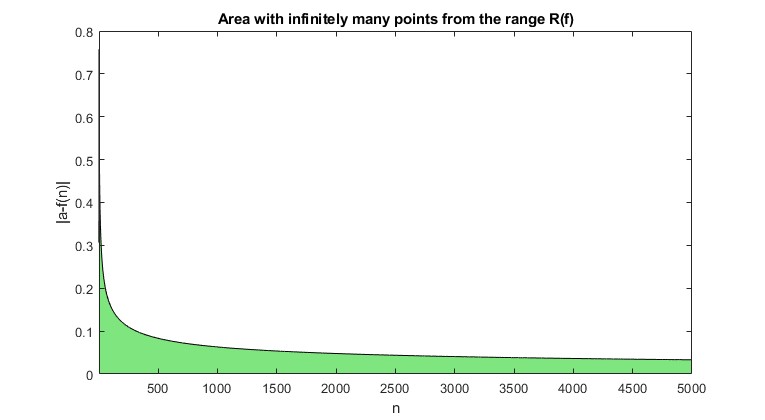}
\begin{center}
 Figure 4: The shadowed area contains infinitely many solutions to 
 inequality (\ref{wolke3}). The narrower this area, the greater the density strength.
\end{center}

\hfill $\qed$

%Let $s=\overline{s}$ be the threshold such that the range of $f_s$, $s<\overline{s}$ is dense and loses the density 
%above $\overline{s}$. The presented computational results below confirm the above analysis and point out that $\overline{s}\approx 1.55$.

%\section{Density strength}
%In \cite{1977Wolke} it is proved that for every $a\geq1$ and every positive $\varepsilon>0$ the inequality
%\[
%\left|f_1(n)-a\right|\leq\frac1{n^{2/5-\varepsilon}}
%\]
%has infinitely many solutions. As an instant corollary we see that the range of $f_1$ is dense in $[1,\infty)$ and the "density strength" is equal $2/5$.
%As shown in the previous section, $f_s$ is dense for $0<s\leq1$ the values of $s$ at which infinite product (\ref{3}) is divergent. 

\section{Density of the complement within rationals}
If $s$ is a natural number, then, clearly, 
the values $f_s(n)$ are rational, $R(f_s)\subset\QQ$. However, there are rational numbers outside of this range.
In fact, if $n$ is composite, then no value $\ds 1+\frac1{m}$ belongs to the range of $f_1$ because 
$\ds f_1(n)=1+\frac1n$ only if $m$ is prime.
% To prove later!
As proved earlier, the range $R(f_1)$ is dense. In this section, we address the question of whether its complement within rationals is also dense. 
Since we are operating within the set of rationals, let us consider only natural values of $s$.

\begin{lm} \label{lm5}
Let $s\in\NN$ and $n$ be composite. Then there is no $m$ such that $\ds f_s(m)=1+\frac1{n^s}$.
\end{lm}
\noi{\bf Proof.} If $m$ were prime, then $\ds f_s(m)=1+\frac1{m^s}=1+\frac1{n^s}$. So, $m=n$ but $n$ is composite. Thus, $m$ is also composite. % Then  $\ds f_s(m) > 1+\frac1{m^s}$, and we have $m>n$. Let $n=\prod q_i$ and $\ds f_s(m)=1+\frac1{n^s}$ for a certain  
Let $m=\prod p_i$. In this product, prime numbers $p_i$ may repeat. 
We have
\[
f_s(m)=\frac{\sum \prod_{i} p_i^s}{m^s}         %\prod_{i=1}^{m} p_i^s}
\]
where the products in the numerator $\prod_{i} p_i^s$ include all combinations of $p_i$, the prime factors of $m$. 
The last expression is equal to $\ds 1+\frac1{n^s}$. Then 
\[
n^s \sum \prod_{i} p_i^s =(n^s+1)m^s.   %\prod_{i=1}^{m} p_i^s.
\]
%At least one addend in the sum $\sum \prod_{j} p_j^s$ does not have
%$p_i^s$ as a multiplier, so $\sum \prod_{j} p_j^s$ is not divisible by $p_i^s$. % \not| p_i$. 
Since $n^s$ and $n^s+1$ are relatively prime, then $n| m$. Therefore the prime factorization of $n$ contains only the primes $p_i$ from the factorization of $m$ and $f_s(n)<f_s(m)$. Therefore 
\[
f_s(m)=1+\frac1{n^s}>f_s(n)>1+\frac1{n^s}.
\]
This contradiction proves Lemma \ref{lm5}. \hfill $\qed$
\begin{cor} \label{cor3}
Let $n$ be a composite number, which is mutually prime with 
$N=\prod\limits_{i=1}^m p_i^{k_i}$. 
Then number %In view of Lemma \ref{lm5}, 
\[
\left(1+\frac1{n^s}\right)f_s(N)
%\cdot\prod_{i=1}^m \left(1+\frac1{p_i^s}\right)
\]
does not belong to the range $R(f_s)$.
%and can be arbitrarily large due to the diverging product as $m\to \infty$.
\end{cor}
% Proof of Corollary 5
%Really, 
%\[
%f_s(N)= \prod\limits_{i=1}^m \left(1+\frac1{p_i^s}+\frac1{p_i^{2s}}+\frac1{p_i^{3s}}+\ldots+
%\frac1{p_i^{sk_i}}\right)=\prod\limits_{i=1}^m \frac{p_i^{s(k_i+1)}-1}{p_i^{s(k_i+1)} - p_i^{ks}},
%\]
%and we apply the argument of Lemma \ref{lm5}.\hfill $\qed$

\begin{thm}
Let $s\geq 0$. Then the rational complement $\QQ\setminus R(f_s)$ is dense in $[1,\infty)$. 
\end{thm}
\noi{\bf Proof.} If $s$ is fractional, then the statement is trivial since $R(f_s)$ is irrational and, consequently, $\QQ\setminus R(f_s)=\QQ$. For $s=0$ the statement is also trivial since $f_0(n)$ is the number of divisors and takes only natural values.
%Let us consider $s=1$. For $(a,b)\subset[1,\infty)$, we choose consecutive primes $\{p_i\}_{i={m_0}}^m$ such that the product \\ 
%$\ds Q_m=\prod_{i=1}^m\left(1+\frac1{p_i}\right)\leq a$ but
%$Q_{m+1}> b$. Then we can find a composite $n>p_{m+1}$ such that $\ds %Q_m\cdot\left(1+\frac1{n}\right)\in(a,b)$. 
%If $Q_m\in(a,b)$, then we select composite $n$ such that
%$\ds Q_m\cdot\left(1+\frac1{n}\right)\in(a,b)$. In any case, we can find a rational number outside the range %$R(f_1)$, which belongs to an arbitrary interval $(a,b)$.

For $s\in\NN$ and an arbitrary interval $(a,b)$ we pick up
a point $x_0\in R(f_s)\cap (a,b)$. Then we choose a sufficiently big composite number $n$, which is mutually prime with $x_0$, such that $\ds x=x_0\cdot \bigl( 1+\frac1{n^s}\bigr) \in (a,b)$. 
In view of Corollary \ref{cor3}, $x\not\in R(f_s)$ but $x\in (a,b)\cap\QQ$.    \hfill $\qed$

\section{Probabilistic estimates: Expectation and Variance}
Let us treat the divisor function $f_s(n)$ as 
a uniformly distributed random variable defined on a sufficiently large interval 
$1\leq n\leq N$. Then, this random variable takes values $f_s(n)$ with probability $\frac1N$.  
Since the values of $n$ are large, then with probability $\frac12$ number $n$
is even,
with probability $\frac13$ number $n$ is divisible by 3, and so on.  
Thus, with probability $\frac1k$ number $n$ has factors $k$ and $\frac{n}{k}$. 
Then 
\[
f_s(n) \approx \frac1{n^s}\left(\frac12 \frac{n^s}{2^s}+ \frac13 \frac{n^s}{3^s} +\ldots
+\frac1k \frac{n^s}{k^s}+\ldots\right).
\]
Then, as $n\to \infty$, 
\bq
f_s(n)\to \sum_{k=1}^\infty \frac{1}{k^{s+1}}=\zeta(s+1).
\eq 
Particularly, the standard divisor function $f_1$ has expectation $\zeta(2)=\frac{\pi^2}{6}$. Thus, we arrived at the following result.
\begin{thm} \label{thmE}
The expectation of the divisor function $f_s$ is equal to $\zeta(s+1)$.
\end{thm}
As a remark, we can observe that as $s=0$, then the expectation becomes infinite, which perfectly fits the sense of $f_0$ as the number of factors. 

For the computational verification of Theorem \ref{thmE}, we use the law of large numbers and calculate the average values of $f_s(n)$,
\begin{equation}
\bar{f}_s=\frac1N \sum_{n=1}^N f_s(n) \to E(f_s)\quad {as}\ N\to\infty.    
\end{equation}
The computational results demonstrate the perfect fit for the presented theory.

%~\hspace{-2.1cm}    \vspace{-1cm} 
\includegraphics[width=156mm, height=9cm]{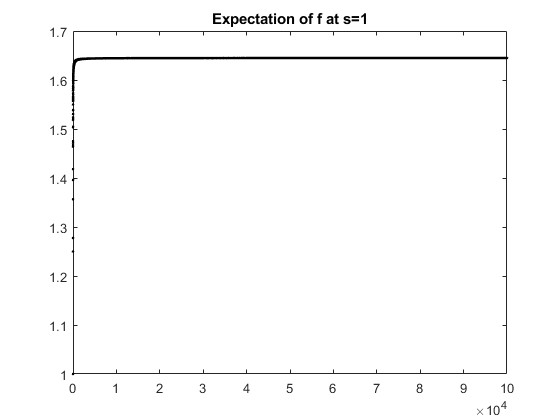}\\[3mm]
\begin{center} Figure 5: Expectation of the divisor function $f_1$ is equal to $\ds \frac{\pi^2}{6}$. 
\end{center}
The dependence of the expectation on $s$ is given in the graph below.\\

\includegraphics[width=156mm, height=9cm]{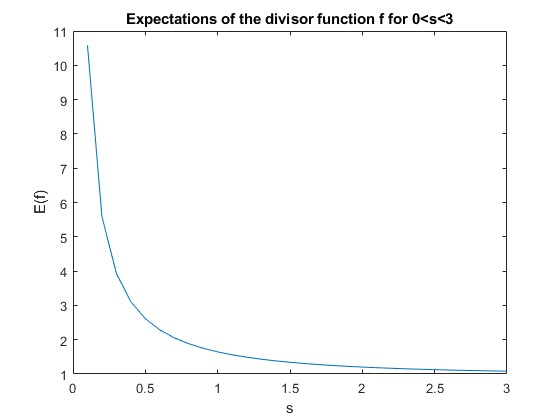}\\[3mm]
\begin{center} Figure 6: Expectation of the divisor function $f_s$ is equal to $\zeta(s+1)$. 
\end{center}

To calculate the variance of random variable $f_s$ for large $n$, we calculate the average of $f_s^2$ and subtract the squared expectation. Then we obtain the decaying to zero function $Var(s)$, which is infinite at $s=0$ and bounded for $s>0$. Particularly, we obtain $Var(1)=0.3$, $Var(2)=0.03$, and $Var(3)=0.005$. Since $f_0(n)$ is just the number of divisors of $n$, and $f_0$ can be drastically changed from 2, if $n=2^k-1$ is prime, to a huge value $k+1$ for $n=2^k$, then it is not surprising that $Var(f_s)\to \infty$ as $s\to0$.

\section{Conclusions}
The paper investigates the dense properties of functions  $\displaystyle f_s(n)=\frac{\sigma_s(n)}{n^s}$ where $\sigma_s$ is the sum of $s$-powers of all factors of $n$. Based on the convergence and divergence properties of the infinite product (\ref{5}), we prove that for every $s>0$, $R(f_s)$, the range of function $f_s$, is locally dense, i.e., any neighborhood of any point from the range contains infinitely many other points from $R(f_s)$. Every number $n$ generates its train $T_s(n)$ -- an ordered collection of ordered sequences $\{f_s(np)\}$.  The proof of local density reveals the structure of the range. It turns out that the range is the union of infinitely many trains. We prove Wolke's conjecture for the points from the range $R_s$ and show that inequality 
\[
|f_s(n)-a|\leq \frac{1}{n^{1-\eps}} 
\]
has infinitely many solutions for any $a\in R_s$ and $\eps>0$. This finding narrows the zone with infinitely many points $f_s$, revealing a stronger quantitative estimate of density strength. 

Then, we address the density for all points in $[1,\infty)$, or the global density. We prove that for $0<s\leq1$, the range of functions $f_s$ is dense and chaotic in $[1,\infty)$ whereas it loses its density properties and experiences ruptures at $s>1$. At threshold $s=1$, $f_s$ loses its chaotic properties. In order to measure density, we extend Wolke's discovery from $s=1$ to all $s\in(0,1]$ and prove the infinite cardinality of the solutions to inequality 
\[
|f_s(n)-a|\leq \frac1{n^{0.4s-\eps}}. 
\]
Our research demonstrates that Wolke's theorem is a particular case of the above inequality.

For integer $s\geq0$, the range is within the set of rationals $\mathbb{Q}$, and for those $s$, we prove that the rational complement to $R(f_s)$, i.e., $\mathbb{Q}\setminus R(f_s)$, is also dense in $[1,\infty)$. For $s=1$, this observation leads to a new partition of the rationals into two dense subsets. 

In the last section, we treat the divisor function as a uniformly distributed random variable and prove that $F(f_s)=\zeta(z+1)$. Particularly, the expectation of the standard divisor function $f_1$ is $\ds \frac{\pi^2}{6}$. It is worth mentioning that $E(f_2)=1.20206$ (Apery's constant) and $\ds E(f_3)=\frac{\pi^4}{90}$. The variance decays when $s$ increases. 
The decay of both the expectation and the variance remains continuous even when the range of the divisor function $f_s$ switches from dense to ruptured at $s=1$.

The analysis is supported by computational demonstrations, highlighting the presented findings.

All of our results show that when considering the dense properties of the divisor function, perfect numbers are not special in any way. While these numbers remain a fascinating topic of study, our research concludes that in this particular context, perfect and multiperfect numbers do not exhibit any distinctive characteristics.

As a remaining open problem, we should point out the question of whether the union of all ranges covers the entire $[1,\infty)$, i.e., whether 
$\bigcup\limits_{s\geq 0} R(f_s)=[1,\infty)$.
Additionally, drawing on the problem of the unknown infiniteness of perfect numbers, we can pose the following question: whether there exists at least one point from the range $R_s(f)$ that repeats infinitely many times. 

Next steps could include the investigation of the density properties of a new function $\displaystyle f_g(n)=\sum\limits_{d|n} \frac{g(d)}{g(n)}$ with increasing function $g$ instead of the power function $\displaystyle f_s(n)=\sum\limits_{d|n} \frac{d^s}{n^s}$. I expect that if $g(n)$ increases slower than the linear function, then the range is also dense. However, the immediate application of the theory developed in this paper is impossible due to the loss of the multiplicative property:  
$f_g(np)\not= f_g(n)\cdot f_g(p)$ unless $g(n)=n^s$. \\\\
{\bf Acknowledgment}

I would like to express my sincere gratitude to Professor L. Boyadzhiev and Mrs. L. Asher for their motivation, support, and guidance.

%\printbibliography[title={Whole bibliography}]

\end{document}